\documentclass[11pt,oneside]{amsart}
\usepackage{fullpage}
\usepackage[left=1in, right=1in]{geometry}
\usepackage{amsmath,amsthm,amsfonts,amssymb,amscd}
\usepackage{relsize}

\makeatletter
\let\save@mathaccent\mathaccent
\newcommand*\if@single[3]{%
	\setbox0\hbox{${\mathaccent"0362{#1}}^H$}%
	\setbox2\hbox{${\mathaccent"0362{\kern0pt#1}}^H$}%
	\ifdim\ht0=\ht2 #3\else #2\fi
}
\newcommand*\rel@kern[1]{\kern#1\dimexpr\macc@kerna}
\newcommand*\widebar[1]{\@ifnextchar^{{\wide@bar{#1}{0}}}{\wide@bar{#1}{1}}}
\newcommand*\wide@bar[2]{\if@single{#1}{\wide@bar@{#1}{#2}{1}}{\wide@bar@{#1}{#2}{2}}}
\newcommand*\wide@bar@[3]{%
	\begingroup
	\def\mathaccent##1##2{%
		\let\mathaccent\save@mathaccent
		\if#32 \let\macc@nucleus\first@char \fi
		\setbox\z@\hbox{$\macc@style{\macc@nucleus}_{}$}%
		\setbox\tw@\hbox{$\macc@style{\macc@nucleus}{}_{}$}%
		\dimen@\wd\tw@
		\advance\dimen@-\wd\z@
		\divide\dimen@ 3
		\@tempdima\wd\tw@
		\advance\@tempdima-\scriptspace
		\divide\@tempdima 10
		\advance\dimen@-\@tempdima
		\ifdim\dimen@>\z@ \dimen@0pt\fi
		\rel@kern{0.6}\kern-\dimen@
		\if#31
		\overline{\rel@kern{-0.6}\kern\dimen@\macc@nucleus\rel@kern{0.4}\kern\dimen@}%
		\advance\dimen@0.4\dimexpr\macc@kerna
		\let\final@kern#2%
		\ifdim\dimen@<\z@ \let\final@kern1\fi
		\if\final@kern1 \kern-\dimen@\fi
		\else
		\overline{\rel@kern{-0.6}\kern\dimen@#1}%
		\fi
	}%
	\macc@depth\@ne
	\let\math@bgroup\@empty \let\math@egroup\macc@set@skewchar
	\mathsurround\z@ \frozen@everymath{\mathgroup\macc@group\relax}%
	\macc@set@skewchar\relax
	\let\mathaccentV\macc@nested@a
	\if#31
	\macc@nested@a\relax111{#1}%
	\else
	\def\gobble@till@marker##1\endmarker{}%
	\futurelet\first@char\gobble@till@marker#1\endmarker
	\ifcat\noexpand\first@char A\else
	\def\first@char{}%
	\fi
	\macc@nested@a\relax111{\first@char}%
	\fi
	\endgroup
}
\makeatother

\makeatletter
\renewcommand*\env@matrix[1][*\c@MaxMatrixCols c]{%
	\hskip -\arraycolsep
	\let\@ifnextchar\new@ifnextchar
	\array{#1}}
\def\thm@space@setup{%
	\thm@preskip=\parskip \thm@postskip=0pt
}
\makeatother

\makeatletter
\@namedef{subjclassname@2020}{\textup{2020} Mathematics Subject Classification}
\makeatother

\usepackage{lastpage}
\usepackage{enumerate}
\usepackage{fancyhdr}
\usepackage{mathrsfs}
\usepackage{xcolor}
\usepackage{graphicx}
\graphicspath{ {./Diagrams/} } 
\usepackage{listings}
\usepackage{hyperref}
\usepackage{euscript}
\usepackage{amsthm}
\usepackage{xifthen}
\usepackage{verbatim}
\usepackage{hyperref}
\usepackage{lipsum}
\usepackage{float}
\usepackage{subfigure}

\usepackage{array}
\newcolumntype{M}[1]{>{\centering\arraybackslash}m{#1}}

\hypersetup{%
	colorlinks=true,
	linkcolor=blue,
	citecolor=blue,
	linkbordercolor={0 0 1}
}

\lstdefinestyle{Python}{
	language        = Python,
	frame           = lines, 
	basicstyle      = \footnotesize,
	keywordstyle    = \color{blue},
	stringstyle     = \color{green},
	commentstyle    = \color{red}\ttfamily
}

\newcommand{\Z}{\mathbb{Z}}
\newcommand{\R}{\mathbb{R}}

\newcommand{\N}{\mathbb{N}}

\newcommand{\sseq}{\subseteq}

\theoremstyle{definition}
\newtheorem{theorem}{Theorem}[section]
\theoremstyle{definition}
\newtheorem{lem}[theorem]{Lemma}
\theoremstyle{definition}

\theoremstyle{definition}
\newtheorem{prop}[theorem]{Proposition}
\theoremstyle{definition}
\newtheorem{cor}[theorem]{Corollary}
\theoremstyle{definition}
\newtheorem{rem}[theorem]{Remark}
\theoremstyle{remark}
\newtheorem{examp}[theorem]{Example}
\theoremstyle{definition}
\newtheorem{defin}[theorem]{Definition}

\author[H.~K.~Chong]{Hip Kuen Chong}
\address{Dept. of Math. \& Stats.\\
	McGill Univ. \\
	Montreal, QC, Canada H3A 0B9 }
\email{chonghk1997@gmail.com}
\author[D.~T.~Wise]{Daniel T. Wise}
\email{wise@math.mcgill.ca}
\date{\today}
\thanks{Research supported by NSERC}

\title{Continuously Many Quasi-isometry Classes of Residually Finite Groups}

\begin{document}

\begin{abstract}
	We study a family of finitely generated residually finite small cancellation groups.  These groups are quotients of $F_2$ depending on a subset $S$ of positive integers.  Varying $S$ yields continuously many groups up to quasi-isometry. 
\end{abstract}

\maketitle

\section{Introduction}

Grigorchuk exhibited continuously many quasiiometry classes of residually finite 3-generator groups by producing continuously many growth types \cite[Thm 7.2]{MR764305}.  
\emph{Continuously many} means having the cardinality of $\R$.  Here, we describe another family of such groups by building upon Bowditch's method for distinguishing quasi-isometry classes \cite{MR1611695}, and use consequences of the theory of special cube complexes to obtain residual finiteness \cite{MR3104553}.  


Consider the rank-$2$ free group $F_2=\langle a,b\rangle$.   Let $w_n = [a,b^{2^{2^n}}][a^2,b^{2^{2^n}}]\cdots [a^{100},b^{2^{2^n}}]$ for  $n\in \N$.

Each subset $S\sseq \N$ is associated to the following group:
$$G(S)\ =\ \langle a,b\mid w_n\colon n\in S\rangle$$ 

In Section~\ref{sec:properties}, we show that $G(S)$ is residually finite when $S\sseq \N_{>100}$.  We also observe that $G(S)$ and $G(S')$ are not quasi-isometric when $S\Delta S'$ is infinite.  


In fact, our proof of residual finiteness for $G(S)$ works in precisely the same way to prove the residual finiteness for the original examples of Bowditch having torsion.  But it appears to fail for Bowditch's torsion-free examples.  We refer to Remark~\ref{rem:BowditchExample}.

We also produced an uncountable family of pairwise non-isomorphic residually finite groups in \cite{ChongWise+2021},   
and perhaps 
an appropriate subfamily also yields continuously many quasi-isometry classes.

Our simple approach arranges for certain infinitely presented small-cancellation groups to be residually finitely presented small-cancellation groups. This approach is likely to permit the construction of other interesting families of finitely generated groups.







\section{Review of Bowditch's Result}

We first recall some small-cancellation background.  See \cite[Ch.V]{smallcancellationbasic}.

\begin{defin}
	For a presentation, a \emph{piece} $p$ is a word appearing in more than one way among the relators.  Note that for a relator $r=q^n$, subwords that differ by a $\Z_n$-action are regarded as appearing in the same way.  A presentation is $C'(\frac{1}{6})$ if $|p|<\frac{1}{6}|r|$ whenever a piece $p$ occurs in a relator $r$.  
\end{defin}

A \emph{major subword} $v$ of a relator $r$ is a subword of a cyclic permutation of $r^{\pm}$ with $|v| > \frac{|r|}{2}$.  Any $g\in G$ is represented by a word $u$ that is \emph{majority-reduced} in the sense that $u$ does not contain a major subword of a relator $r$.  We will use the following proposition for $C'(\frac{1}{6})$ groups, which are groups given by a $C'(\frac{1}{6})$ presentation.  

\begin{prop}\label{prop:dehn}
	Let $G$ be a $C'(\frac{1}{6})$ group.  If $u$ is a nonempty majority-reduced word, then $u\neq 1_G$.
\end{prop}

We now recall definitions leading to the main theorem of \cite{MR1611695}. 
Let $\N^+ = \{n\in \Z\colon n\geq 1\}$.


\begin{defin}\label{def:krelated}
	Two subsets $L, L'\sseq \N^+$ are \emph{related} if for some $k\geq 1$:
	\begin{enumerate}
		\item for any $m\in L$ with $m>(k+1)^2$, there is $m'\in L'$ with $m'\in [\frac{m}{k}, km]$; and
		\item for any $m'\in L'$ with $m'>(k+1)^2$, there is $m\in L$ with $m\in [\frac{m'}{k}, km']$.
	\end{enumerate}
	
	We write $L\sim L'$ if $L$ and $L'$ are related, and write $L\not\sim L'$ otherwise.
\end{defin}


\begin{lem}
	The relation $\sim$ in Definition~\ref{def:krelated} is an equivalence relation on subsets of $\N^+$.  
\end{lem}

\begin{proof}
	The relation $\sim$ is reflexive by letting $k=1$.  The relation $\sim$ is symmetric by definition.  Hence it suffices to show $\sim$ is transitive.
	
	Let $S, S', S''\sseq \N^+$.  Suppose $S\sim S'$ via $k$ and $S'\sim S''$ via $k'$.  We claim that $S\sim S''$ via $kk'$.  Let $m\in S$.  There is $m'\in S'$ with $m'\in [\frac{m}{k},km]$ and there is $m''\in S''$ with $m''\in [\frac{m'}{k'},k'm']$.  Thus, $m''\in \left[\frac{m}{kk'}, kk'm\right]$. 
	Similarly, $m\in \left[\frac{m''}{kk'}, kk'm''\right]$ for $m''\in S''$.
\end{proof}

For sets $S,S'$, their \emph{symmetric difference}  is  $S\Delta S' = (S - S')\cup (S' - S)$.

\begin{examp}\label{examp:Bowditch22n}
	If $S,S'\sseq \N^+$ with infinite $S\Delta S'$, then $\{2^{2^n}\}_{n\in S} \not\sim \{2^{2^m}\}_{m\in S'}$  {\cite[Lem 4]{MR1611695}}.
\end{examp}


With the notion of $\sim$, the following is a simplified version of the main theorem in \cite{MR1611695}.

\begin{theorem}\label{thm:main}
	Let $G$ and $G'$ be the finitely generated $C'(\frac{1}{6})$ groups presentated below.\\  If $G$ is quasi-isometric to $G'$, then $\{|w_i|\}_{i\in I} \sim \{|w_j'|\}_{j\in J}$.  
	$$G\ =\ \langle A\mid w_i\colon i\in I\rangle, \qquad G'\ =\ \langle A\mid w_j'\colon j\in J\rangle.$$
\end{theorem}

\section{Proving the Family of Groups Have Desired Properties}\label{sec:properties}

\subsection{Small Cancellation}

Let $\N_{>100}=\{n\in \N\colon n>100\}$.

\begin{prop}\label{prop:smallcan}
	For any infinite subset $S\sseq \N_{>100}$, the associated group $G$ is $C'(\frac{1}{6})$.  \\
	Furthermore, $G_k\ =\ \left\langle a,b\ \middle|\ b^{2^{2^{k}}}, w_n\colon n\in S, n<k\right\rangle$ is $C'(\frac{1}{6})$ for each $k\in \N$. 
\end{prop}

\begin{proof}
	For the first statement, it suffices to show that $w_n$ and $w_m$ have small overlap for $n>m>100$.  The longest piece between $w_n$ and $ w_m$ is $b^{-2^{2^{m}}}a^{100}b^{2^{2^{m}}}$.  Thus, $C'(\frac{1}{6})$ holds since:
	$$\left|b^{-2^{2^{m}}}a^{100}b^{2^{2^{m}}}\right| 
	\ =\ 100+2\cdot 2^{2^{m}}
	\ <\ \frac{1}{6}\left(10100+200\cdot 2^{2^{m}}\right) 
	\ =\ \frac{1}{6}|w_m|
	\ <\ \frac{1}{6}|w_n|$$
	
	For the second statement, we additionally show that $w_n$ and $b^{2^{2^k}}$ satisfy the $C'(\frac{1}{6})$ condition for $100<n<k$.  Their longest piece is $b^{2^{2^n}}$, which is shorter than $\frac{1}{6}$ of the lengths of $w_n$ and $b^{2^{2^k}}$.
\end{proof}

\subsection{Residual Finiteness}

Observe that $G_k=G/\langle\!\langle b^{2^{2^k}}\rangle\!\rangle$ since $w_m\in \langle\!\langle b^{2^{2^k}}\rangle\!\rangle$ for $m\geq k$.  

\begin{prop}\label{prop:ResResSmallCan}
	For any infinite subset $S\sseq \N_{>100}$, the associated group $G$ is residually finite.
\end{prop}

\begin{proof}
 Since $G_k$ is a finitely presented $C'(\frac{1}{6})$ group, the hyperbolic group $G_k$ is cocompactly cubulated by \cite{MR2053602}. Thus $G_k$ is residually finite by \cite{MR3104553}.  
	
	By Proposition~\ref{prop:dehn}, each $g\in G-\{1\}$ is represented by a majority-reduced word $v$, which does not contain a majority subword of $b^{2^{2^{|v|}}}$ since $|v|<\frac{1}{2}\cdot 2^{2^{|v|}}=\frac{1}{2}\left|b^{2^{2^{|v|}}}\right|$.  Hence, $v\neq 1_{G_{|v|}}$.
	Thus, $G$ is residually residually finite and hence residually finite.  
\end{proof}

\begin{rem}\label{rem:BowditchExample}
	Bowditch's original examples were $B(S) = \langle a, b \mid \big(a^{2^{2^n}}b^{2^{2^n}}\big)^7 \colon n\in S\sseq \N\rangle$.  As in Proposition~\ref{prop:ResResSmallCan}, $B(S)$ is residually finite since it is residually finitely presented $C'(\frac{1}{6})$ using the quotients to $B/\langle\!\langle a^{2^{2^n}}, b^{2^{2^n}}\rangle\!\rangle$ for $n\geq 3$.  However, the analogous argument fails for Bowditch's torsion-free examples $B'(S) = \langle a, b \mid a\big(a^{2^{2^n}}b^{2^{2^n}}\big)^{12} \colon n\in S\sseq \N\rangle$.  
\end{rem}

\subsection{Pairwise Non-quasi-isometric}
We first prove a lemma about the relation $\sim$.  


\begin{lem}\label{lem:krelated}
	$S\sim nS\sim (S+n)$ for $n\in \N^+$ and $S\sseq \N^+$.
\end{lem}
\begin{proof}
	First, $S\sim nS$ via $n$.  Indeed, for any $s\in S$, $ns\in \left[\frac{s}{n},ns\right]$; for any $ns\in nS$, $s\in \left[\frac{ns}{n},n\cdot ns\right]$. 
	
	Moreover, $S\sim (S+n)$ via $n+1$.  For any $s\in S$, $s+n\leq (n+1)s$, so $s+n\in \left[\frac{s}{n+1},(n+1)s\right]$.  On the other hand, for $s+n\in S+n$, $(n+1)s\geq s+n$ implies $s\geq \frac{s+n}{n+1}$.  Hence, $s\in \left[\frac{s+n}{n+1}, (n+1)(s+n)\right]$.  
\end{proof}

\begin{prop}
	Let $S,S'\sseq \N^+$ have infinite $S\Delta S'$, then $\{|w_n|\}_{n\in S}\not\sim \{|w_m|\}_{m\in S'}$.
\end{prop}

\begin{proof}
	$\{|w_n|\colon n\in S\} = \{10100+200\cdot 2^{2^n}\colon n\in S\} = 10100+200\cdot \{2^{2^n}\colon n\in S\}$.  By Lemma~\ref{lem:krelated}, $\{|w_n|\colon n\in S\}\sim \{2^{2^n}\colon n\in S\}$. Similarly, $\{|w_m|\colon m\in S'\}\sim \{2^{2^m}\colon m\in S'\}$.  By Example~\ref{examp:Bowditch22n}, $\{2^{2^n}\colon n\in S\}\not\sim \{2^{2^m}\colon m\in S'\}$,  so $\{|w_n|\}_{n\in S}\not\sim \{|w_m|\}_{m\in S'}$.  
\end{proof}

\begin{cor}
	If $S,S'\sseq \N_{>100}$ have infinite $S\Delta S'$, then $G(S)$ and $G(S')$ are not quasi-isometric.
\end{cor}
\begin{proof}
	$\{|w_n|\}_{n\in S}\not\sim \{|w_m|\}_{m\in S'}$, hence $G(S)$ and $G(S')$ are not quasi-isometric by Theorem~\ref{thm:main}.
\end{proof}


For $A,B\sseq N$, declare $A\sim_{_\Delta}\!B$ if $\left|A\Delta B\right|<\infty$.  As noted by Bowditch, each $\sim_{_\Delta}\!$ equivalence class is countable.  Hence, there are continuously many $\sim_{_\Delta}\!$ equivalence classes.  Our construction thus produces continuously many pairwise non-quasi-isometric groups $G(S)$, which are $C'(\frac{1}{6})$ and residually finite.

\bibliographystyle{alpha}
\bibliography{Citation}
\end{document}